\documentclass[12pt,a4paper]{article}
\pdfoutput=1

\usepackage[T1]{fontenc}
\usepackage[utf8]{inputenc}
\usepackage[UKenglish]{isodate}
\cleanlookdateon

\usepackage{amsfonts,amsmath,amssymb,amsthm,mathtools,dsfont,microtype} 
\usepackage[usenames,dvipsnames]{xcolor} 
\usepackage[noBBpl]{mathpazo} 

\usepackage[labelsep=period,labelfont=bf,justification=centering]{caption}
\usepackage{float,graphicx,subcaption}
\usepackage{booktabs,makecell}
\usepackage{enumitem,mdwlist}
\setlist[itemize]{topsep=0ex,itemsep=0ex,parsep=0.4ex}
\setlist[enumerate]{topsep=0ex,itemsep=0ex,parsep=0.4ex}

\usepackage{parskip,fullpage}
\usepackage{standalone}
\usepackage{thm-restate}

\usepackage[hyphens]{url} 
\usepackage[linktoc=all,hidelinks,colorlinks,unicode=true]{hyperref} 
\usepackage[capitalise,compress,nameinlink,noabbrev]{cleveref} 
\hypersetup{linkcolor={blue!70!black},citecolor={black},urlcolor={blue!70!black}}

\usepackage{tikz}

\newtheorem{theorem}{Theorem}
\newtheorem{lemma}[theorem]{Lemma}

\newtheorem*{claim*}{Claim}
\newtheorem{conjecture}[theorem]{Conjecture}

\theoremstyle{definition}
\newtheorem{remark}[theorem]{Remark}
\newtheorem{problem}[theorem]{Problem}

\newenvironment{poc}{\begin{proof}}{\end{proof}}

\newcommand{\defn}[1]{\textcolor{Maroon}{\emph{#1}}}

\newcommand*{\bN}{\mathbb{N}}

\newcommand*{\cF}{\mathcal{F}}

\newcommand*{\cJ}{\mathcal{J}}

\newcommand*{\cC}{\mathcal{C}}

\newcommand*{\cG}{\mathcal{G}}

\DeclareMathOperator{\sq}{\square}

\title{\texorpdfstring{\vspace{-4ex}}{}A note on graphs of $k$-colourings}

\date{\today}

\author{
Emma Hogan\footnotemark[2],\ \quad Alex Scott\footnotemark[2]\ \footnotemark[1], \quad
Youri Tamitegama\footnotemark[2],\ \quad
Jane Tan\footnotemark[2]\ \footnotemark[3]
}

\begin{document}
\maketitle

\renewcommand{\thefootnote}{\fnsymbol{footnote}} 
\footnotetext[1]{Supported by EPSRC grant EP/X013642/1.}
\footnotetext[2]{Mathematical Institute, University of Oxford, United Kingdom
(\textsf{\{\href{mailto:hogan@maths.ox.ac.uk}{hogan},\href{mailto:scott@maths.ox.ac.uk}{scott},\href{mailto:tamitegama@maths.ox.ac.uk}{tamitegama},\href{mailto:jane.tan@maths.ox.ac.uk}{jane.tan}\}\allowbreak @maths.ox.ac.\allowbreak uk}).}
\footnotetext[3]{All Souls College, University of Oxford.}

\renewcommand{\thefootnote}{\arabic{footnote}} 
\begin{abstract}

For a graph $G$, the $k$-colouring graph of $G$ has vertices corresponding to proper $k$-colourings of $G$ and edges between colourings that differ at a single vertex. 
The graph supports the Glauber dynamics Markov chain for $k$-colourings, and has been extensively studied from both extremal and probabilistic perspectives. 

In this note, we show that for every graph $G$, there exists $k$ such that $G$ is uniquely determined by its $k$-colouring graph, confirming two conjectures of Asgarli, Krehbiel, Levinson and Russell. We further show that no finite family of generalised chromatic polynomials for $G$, which encode induced subgraph counts of its colouring graphs, uniquely determine $G$.
\end{abstract}

\section{Introduction}
Let $G$ be a graph on vertex set $V(G)$ and edge set $E(G)$. 
Throughout this paper, all colourings are proper, and a \defn{$k$-colouring} is a proper colouring using at most $k$ colours from a fixed palette, say $[k] \coloneqq \{1,\dotsc,k\}$. 
The \defn{chromatic polynomial} $\pi_G(k)$ counts the number of $k$-colourings of $G$ as a function of $k$. 
Chromatic polynomials were first considered for planar maps by Birkhoff~\cite{B12} in 1912, and then for arbitrary graphs by Whitney~\cite{W32} in 1932. Since then, they have been well-studied in the literature, with considerable interest in ways in which they can be computed, their algebraic properties, and generalisations (see~\cite{Read68} for a classical introduction). 

A more detailed picture of the set of $k$-colourings of a graph $G$ is given by the \defn{$k$-colouring graph} $\cC_k(G)$: this has vertex set the $k$-colourings of $G$, and edges between pairs of $k$-colourings that differ at precisely one vertex of $G$.
Random walks on the $k$-colouring graph give the Glauber dynamics Markov chain,
which has been extensively studied from the perspective of random sampling and approximate counting of $k$-colourings (see for example \cite{DFFE06, Jerrum95, Vigoda00}).  The $k$-colouring graph has also been investigated in the context of combinatorial reconfiguration
(see, for example, the surveys in~\cite[Chapter 10]{50years} and \cite{heuvel2013complexity}). 

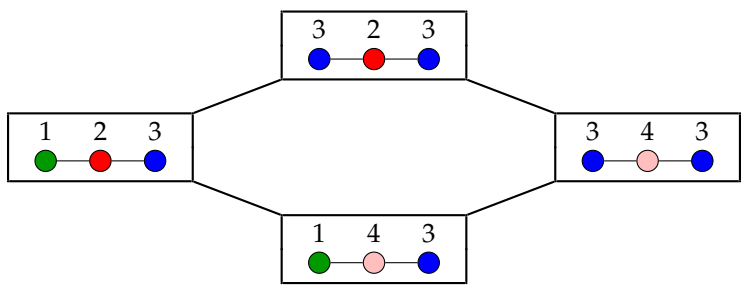
\begin{figure}[ht]
 \footnotesize    
     \centering
     \begin{tikzpicture}[
    scale=0.9,
    mycircle/.style={
        circle,
        draw=black,
        fill=black,
        fill opacity = 1,
        inner sep=0pt,
        minimum size=8pt,
        font=\small},
    nocircle/.style={
        circle,
        draw=black,
        fill=black,
        fill opacity = 1,
        inner sep=0pt,
        minimum size=0.45pt,
        font=\small},
    myarrow/.style={-},
    dottedarrow/.style={-,dashed},
    thiccarrow/.style={-,line width=0.9pt},
    node distance=1.2cm and 1.5cm
]

\begin{scope}
    \begin{scope}[]
        \foreach \x/\y/\name in {-1.35/0/a, 1.35/0/b1, 1.35/0.5/x, 1.35/1/c1, -1.35/1/d}{ 
        \node[nocircle] (\name) at (\x,\y){};
        }
    \end{scope}

    \node[mycircle, fill=green!60!black, label={$1$}] (1) at (-0.8,0.3){};
    \node[mycircle, fill=red, label={$2$}] (2) at (0,0.3){};
    \node[mycircle, fill=blue, label={$3$}] (3) at (0.8,0.3){};

    \path[every node/.style={font=\sffamily\small}]
        (1) edge [color=black] (2)
        (2) edge [color=black] (3)
        (a) edge [color=black, thick ] (b1)
        (b1) edge [color=black, thick ] (c1)
        (c1) edge [color=black, thick ] (d)
        (d) edge [color=black, thick ] (a);
\end{scope}

\begin{scope}[xshift=4cm, yshift=-1.5cm]
    \begin{scope}[]
        \foreach \x/\y/\name in {-1.35/0/a, -1.35/0.5/y, 1.35/0.5/y2, 1.35/0/b, 1.35/1/c2, -1.35/1/d2}{ 
        \node[nocircle] (\name) at (\x,\y){};
        }
    \end{scope}

    \node[mycircle, fill=green!60!black, label={$1$}] (1) at (-0.8,0.3){};
    \node[mycircle, fill=pink, label={$4$}] (2) at (0,0.3){};
    \node[mycircle, fill=blue, label={$3$}] (3) at (0.8,0.3){};

    \path[every node/.style={font=\sffamily\small}]
        (1) edge [color=black] (2)
        (2) edge [color=black] (3)
        (a) edge [color=black, thick ] (b)
        (b) edge [color=black, thick ] (c2)
        (c2) edge [color=black, thick ] (d2)
        (d2) edge [color=black, thick ] (a);
\end{scope}

\begin{scope}[xshift=4cm, yshift=1.5cm]
    \begin{scope}[]
        \foreach \x/\y/\name in {-1.35/0/a3, -1.35/0.5/z, 1.35/0.5/z2, 1.35/0/b3, 1.35/1/c, -1.35/1/d}{ 
        \node[nocircle] (\name) at (\x,\y){};
        }
    \end{scope}

    \node[mycircle, fill=blue, label={$3$}] (1) at (-0.8,0.3){};
    \node[mycircle, fill=red, label={$2$}] (2) at (0,0.3){};
    \node[mycircle, fill=blue, label={$3$}] (3) at (0.8,0.3){};

    \path[every node/.style={font=\sffamily\small}]
        (1) edge [color=black] (2)
        (2) edge [color=black] (3)
        (a3) edge [color=black, thick ] (b3)
        (b3) edge [color=black, thick ] (c)
        (c) edge [color=black, thick ] (d)
        (d) edge [color=black, thick ] (a3);
\end{scope}

\begin{scope}[xshift=8cm]
    \begin{scope}[]
        \foreach \x/\y/\name in {-1.35/0/a4, -1.35/0.5/w, 1.35/0/b, 1.35/1/c, -1.35/1/d4}{ 
        \node[nocircle] (\name) at (\x,\y){};
        }
    \end{scope}

    \node[mycircle, fill=blue, label={$3$}] (1) at (-0.8,0.3){};
    \node[mycircle, fill=pink, label={$4$}] (2) at (0,0.3){};
    \node[mycircle, fill=blue, label={$3$}] (3) at (0.8,0.3){};

    \path[every node/.style={font=\sffamily\small}]
        (1) edge [color=black] (2)
        (2) edge [color=black] (3)
        (a4) edge [color=black, thick ] (b)
        (b) edge [color=black, thick ] (c)
        (c) edge [color=black, thick ] (d4)
        (d4) edge [color=black, thick ] (a4);
\end{scope}


    \path[every node/.style={font=\sffamily\small}]
    (b1) edge [color=black, thick ] (d2)
    (c1) edge [color=black, thick ] (a3)
    (c2) edge [color=black, thick ] (a4)
    (b3) edge [color=black, thick ] (d4);

\end{tikzpicture}
     \vspace{-0.5cm}
     \caption{An induced $C_4$ in $\cC_k(P_3)$ for $k\geq 4$.}
     \label{fig:c4p3}
\end{figure}

The chromatic polynomial $\pi_G(k)$ counts the number of vertices in $\cC_k(G)$.
Asgarli, Krehbiel, Levinson and Russell~\cite{AKLR24} recently introduced a more general family of functions by replacing vertex counts with counts of instances of a fixed arbitrary graph: for graphs $G$ and $H$, and $k\in\bN$, the
\defn{generalised chromatic polynomial} $\pi^{(H)}_G(k)$ is 
the number of subsets of $V(\cC_k(G))$ that induce a subgraph isomorphic to $H$ as a function of $k$. Thus $\pi^{(K_1)}_G(k)$ is the chromatic polynomial of $G$, and $\pi^{(K_2)}_G(k)$ counts the number of edges in $\cC_k(G)$ (see~\cref{fig:c4p3} for an example of a structure that contributes to $\pi^{(C_4)}_{P_3}(4)$). 

For fixed graphs $G$ and $H$,
Asgarli et al.~proved that  $\pi^{(H)}_G(k)$ is a polynomial in $k$ sufficiently large relative to the size of $H$. In \cref{sec:poly}, we strengthen this result to show that $\pi^{(H)}_G(k)$ is a polynomial without restriction.
\begin{restatable}{theorem}{thmpoly}\label{thm:poly}
For any fixed graphs $G$ and $H$, the function
$\pi^{(H)}_G(k)$ is a polynomial in $k$.
\end{restatable}

Asgarli et al. also discuss the extent to which a graph $G$ is determined by the invariants~$\pi_G^{(H)}$. Letting $\cG$ be the set of finite graphs, they conjecture that the collection of polynomials $\{\pi^{(H)}_G(k)\}_{H\in \cG}$, or equivalently the collection of all colouring graphs $\{\cC_k(G)\}_{k\in\bN}$, is a complete graph invariant.
\begin{conjecture}[Conjecture 6.1 \cite{AKLR24}] \label{conj:2}
    For any graph $G$, the collection $\{\cC_k(G)\}_{k\in\bN}$ uniquely determines $G$.
\end{conjecture}

They also make the stronger conjecture that finitely many colouring graphs suffice.
\begin{conjecture}[Conjecture 6.2 \cite{AKLR24}]
\label{conj:3}
    There exists some function $f: \cG \to \bN$ such that for any graph $G$, the collection $\{\cC_k(G)\}^{f(G)}_{k=1}$ uniquely determines $G$.
\end{conjecture}
In \cref{sec:inv}, we confirm both conjectures by proving a stronger result.
\begin{restatable}{theorem}{reconstructingG}
\label{thm:reconstructing_G}
     Let $G$ be a graph on $n$ vertices. For any natural number $k > 5n^2$, the colouring graph $\cC_k(G)$ uniquely determines $G$ up to isomorphism.
\end{restatable}

Since the collection of all colouring graphs $\{\cC_k(G)\}_{k\in\bN}$ holds the same information as the collection of generalised chromatic polynomials $\{\pi^{(H)}_G(k)\}_{H\in \cG}$, another natural direction to investigate is whether a finite subcollection of generalised chromatic polynomials suffices to distinguish all non-isomorphic graphs.
In \cref{sec:finite} we give a negative answer. 

\begin{restatable}{theorem}{finitefamily}
\label{thm:no_finite_family}
    No finite family of generalised chromatic polynomials is a complete graph invariant.
\end{restatable}

\section{Polynomiality}
\label{sec:poly}

Let $\pi^{(H)}_G(k)$ denote the number of induced copies of $H$ in $\cC_k(G)$. We extend the standard proof of polynomiality for chromatic polynomials via partitions to generalised chromatic polynomials.

\thmpoly*

\begin{proof}
Let $h=|H|$ and $n=|G|$.
We will say that a partition $P_1\cup\dots\cup P_s$ of $V(G)\times [h]$ is \defn{valid} if $P_j\cap (V(G)\times \{i\})$ is an independent set in $G$ for each $i\in [h]$ and each $j$.
Fix an ordering~$\prec$ on the vertices of $\cC_k(G)$.
Any collection $S = \{ c_1\prec \dotsc \prec c_h\}$ of $k$-colourings of $G$ defines a function $c\colon V(G)\times [h]\rightarrow [k]$ by $c(v,i) = c_i(v)$ for each $v\in V(G)$, $i\in [h]$, and $c$ induces a valid partition $P_c = \{ c^{-1}(i) \colon i\in [k]\}$ of $V(G)\times [h]$.

The graph induced by $S$ in $\cC_k(G)$ depends only on $P_c$, in the sense that if another collection of $h$ colourings $S'$ defines a partition $P_{c'}$ of $V(G)\times [h]$ then $S$ and $S'$ induce isomorphic subgraphs of $\cC_k(G)$ if $P_c=P_{c'}$.
Each partition $P$ of $V(G)\times [h]$ thus corresponds to a fixed induced graph.
Provided $P$ is valid and consists of $t$ non-empty parts, we can colour its parts in $(k)_t$ different ways, where $(k)_t \coloneqq k(k-1)\dotsb(k-t+1)$ denotes the falling factorial. 
Thus, each such partition that yields an induced copy of $H$ contributes exactly $(k)_t$ induced copies of $H$ to $\cC_k(G)$.
Writing $N_t^{(H)}$ for the number of valid partitions 
with exactly $t$ parts
yielding an induced copy of $H$, the generalised chromatic polynomial of $H$ is given by the formula
\[
\pi^{(H)}_G(k) = \sum_{t=1}^{n} N_t^{(H)} {\binom{k}{t}} t! = \sum_{t=1}^{n} N_t^{(H)} (k)_t\ .
\]
Since each summand is a polynomial and $n$ is fixed, $\pi_G^{(H)}$ is a polynomial as well.
\end{proof}

\newpage
\section{Complete invariance}
\label{sec:inv}
We now prove that the collection of colouring graphs gives a complete graph invariant.
Let $G$ be a graph on $n$ vertices.
A vertex $c\in V(\cC_k(G))$ is \defn{rainbow} if it represents a colouring of $G$ using $n$ distinct colours; that is, $c(u) \neq c(v)$ for any distinct $u$ and $v$ in $V(G)$. 
Our strategy for reconstructing $G$ is to choose a $\cC_k(G)$ with $k$ large enough so that 
most vertices of $C_k(G)$ correspond to a rainbow colouring; we will then be able to use the clique structure to reconstruct the graph.

\begin{lemma}\label{lem:triangles}
Let $G$ be a graph on $n$ vertices.
If $c_1,c_2,c_3$ are vertices of $\cC_k(G)$ inducing a copy of $K_3$, then $c_1,c_2,c_3$ differ as colourings at a single vertex of $G$.
\end{lemma}
\begin{proof}
Suppose that $c_1$ and $c_2$ differ at vertex $u$ of $G$, while $c_2$ and $c_3$ differ at vertex $v$ with $u\neq v$.
Then $c_1$ and $c_3$ differ at both vertices $u$ and $v$, so $c_1c_3$ is not an edge of $\cC_k(G)$, a contradiction.
\end{proof}
It follows that vertices in any clique in $\cC_k(G)$ correspond to colourings which differ at a single vertex $v$ of $G$.
We say that such cliques are \defn{generated} by $v$.
For a colouring $c$ in $\cC_k(G)$, let $\cJ(c)$ be the collection of maximal cliques containing $c$ in $\cC_k(G)$.
When $k\geq n+3$, $\cJ(c)$ consists of $n$ cliques, each generated by a distinct vertex of $G$.
Say that $c$ is \defn{typical} if for each $v\in G$, the clique generated by $v$ in $\cJ(c)$ is of size $k-\deg(v)$.
We note that every rainbow colouring is typical.

\begin{lemma}\label{lem:3n2_bound}
Let $G$ be a graph on $n$ vertices and take any natural number $k>3n^2$. Then $\cC_k(G)$ uniquely determines the degree sequence of $G$.
Moreover, more than half of all vertices in $\cC_k(G)$ are rainbow.
\end{lemma}
\begin{proof}
The number of $k$-colourings of $G$ is at most $k^n$, and there are $\binom{k}{n}\cdot n!\geq k^n(1-\frac{n}{k})^n$ colourings of $G$ using $n$ colours.
Since $k>3n^2$, the proportion of vertices in $\cC_k(G)$ that are rainbow is at least $(1-\frac{n}{k})^n > 1/2$. 
For each vertex $c\in \cC_k(G)$, 
we consider the sizes of maximal cliques containing $c$.  The majority will be typical and therefore give the same collection of sizes.
From any such typical vertex $c$, we can then deduce the degree sequence of $G$ by subtracting the size of each maximal clique containing $c$ from $k$.
\end{proof}

\reconstructingG*

\begin{proof}
    Consider a vertex $c$ of $\cC_k(G)$, and let $J_u, J_v\in \cJ(c)$ be two cliques generated by distinct vertices $u,v\in V(G)$ respectively. 
    We will show that, by counting $4$-cycles that contain $c$ and intersect $J_u\setminus c$ and $J_v\setminus c$, we can determine whether $u$ and $v$ are adjacent in~$G$ (see \cref{fig:maximal_cliques}).
    \newpage
    \begin{figure}[ht]
     \centering
     \begin{tikzpicture}[
    mycircle/.style={
        circle,
        draw=black,
        fill=black,
        fill opacity = 1,
        inner sep=0pt,
        minimum size=5pt,
        font=\small},
    smallcircle/.style={
        circle,
        draw=black,
        fill=black,
        fill opacity = 1,
        inner sep=0pt,
        minimum size=3pt,
        font=\small},
    nocircle/.style={
        circle,
        draw=black,
        fill=black,
        fill opacity = 0,
        inner sep=0pt,
        minimum size=0pt,
        font=\small},
    myarrow/.style={-},
    dottedarrow/.style={-,dashed},
    thiccarrow/.style={-,line width=0.9pt},
    node distance=1.2cm and 1.5cm
]

\tikzstyle{every node}=[font=\small]

\begin{scope}
    \node[mycircle, label=left:{$c_0$}] (1) at (-2,0){};
    \foreach \x/\y/\name in {0/2.5/a1,0/1.5/a3, 0/0.7/b1, 0/-0.3/b3, 0/-1.5/c1, 0/-2.5/c3}{ 
        \node[nocircle] (\name) at (\x,\y){};
    }

    \foreach \x/\y/\name in {1/2.4/va1, 1/2.2/va2, 1/2/va3, 1/1.6/va4, 1/0.5/vb1, 1/0.3/vb2, 1/-0.1/vb3, 1/-1.7/vc1, 1/-2.1/vc2, 1/-2.3/vc3}{ 
        \node[smallcircle] (\name) at (\x,\y){};
        \path[every node/.style={font=\sffamily\small}]
        (1) edge [color=black] (\name);
    }

    \node[rotate=90, scale=1.4, blue] () at (1,-1.1) {...};
    \node[rotate=90, scale=0.8] () at (1,0.1) {...};
    \node[rotate=90, scale=0.8] () at (1,-1.9) {...};
    \node[rotate=90, scale=0.8] () at (1,1.8) {...};

    \foreach \x/\y/\name in {2/2.5/x1,  2/1.5/x3}{ 
        \node[nocircle] (\name) at (\x,\y){};
    }
    \path[every node/.style={font=\sffamily\small}]
        (a1) edge [color=black, dashed, blue] (a3)
        (a3) edge [color=black, dashed, blue] (x3)
        (x3) edge [color=black, dashed, blue] (x1)
        (x1) edge [color=black, dashed, blue] (a1);

    \foreach \x/\y/\name in {2/0.7/y1,  2/-0.3/y3}{ 
        \node[nocircle] (\name) at (\x,\y){};
    }
    \path[every node/.style={font=\sffamily\small}]
        (b1) edge [color=black, dashed, blue] (b3)
        (b3) edge [color=black, dashed, blue] (y3)
        (y3) edge [color=black, dashed, blue] (y1)
        (y1) edge [color=black, dashed, blue] (b1);

    \foreach \x/\y/\name in {2/-1.5/z1,  2/-2.5/z3}{ 
        \node[nocircle] (\name) at (\x,\y){};
    }
    \path[every node/.style={font=\sffamily\small}]
        (c1) edge [color=black, dashed, blue] (c3)
        (c3) edge [color=black, dashed, blue] (z3)
        (z3) edge [color=black, dashed, blue] (z1)
        (z1) edge [color=black, dashed, blue] (c1);

    \node[smallcircle] (c21) at (3,1.5){};
    \node[smallcircle] (c22) at (3,1.3){};
    \node[smallcircle] (c23) at (3,1.1){};
    \node[smallcircle] (c24) at (3,0.7){};
    \node[smallcircle] (c25) at (3,-1){};

    \path[every node/.style={font=\sffamily\small}]
        (c21) edge [color=black] (va1)
        (c21) edge [color=black] (vb1)
        (c22) edge [color=black] (va2)
        (c22) edge [color=black] (vb3)
        (c23) edge [color=black] (va3)
        (c23) edge [color=black] (vb2)
        (c24) edge [color=black] (va4)
        (c24) edge [color=black] (vb3)
        (c25) edge [color=black] (va4)
        (c25) edge [color=black] (vc1);

    \node[smallcircle, label={[left, blue]$v_1$}, blue] (v1) at (5.5,2){};
    \node[smallcircle, label={[left, blue]$v_2$}, blue] (v2) at (5.5,0.2){};
    \node[smallcircle, label={[left, blue]$v_n$}, blue] (v3) at (5.5,-2){};
    \node[rotate=90, scale=1.2, blue] () at (5.5,-0.8) {$\dotsc$};
    
    \path[every node/.style={font=\sffamily\small}]
        (v1) edge [color=black, bend left=20, blue] (v2)
        (v1) edge [color=black, dashed, bend left=25, blue] (v3);

    \draw[blue] (1,2.75) node {$J_{v_1} \setminus c_0$};
    \draw[blue] (1,-0.55) node {$J_{v_2} \setminus c_0$};
    \draw[blue] (1,-2.75) node {$J_{v_n} \setminus c_0$};

    \draw (0.5,-3.5) node {$\mathcal{C}_k(G)$};
    \draw (5.5,-3.5) node {$G$};
\end{scope}

\end{tikzpicture}
     \vspace{-1cm}
     \caption{Detecting edges in $G$ by counting $4$-cycles containing $c_0$ in $\cC_k(G)$.}
     \label{fig:maximal_cliques}
 \end{figure}
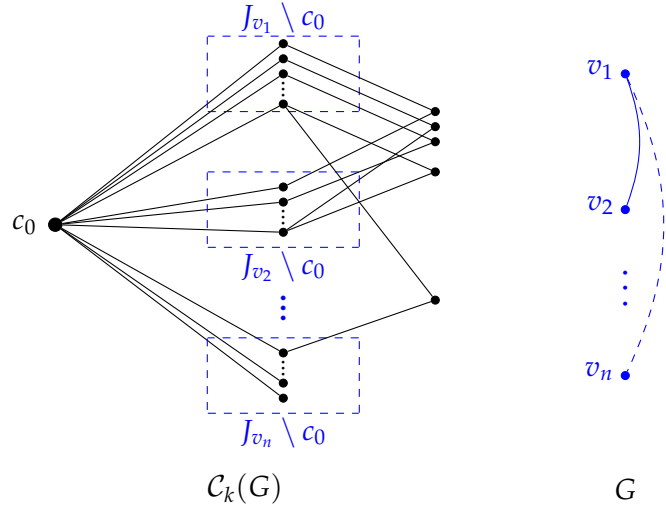
    \begin{claim*}
        Let $c_0$ be rainbow and $J_u,J_v\in \cJ(c_0)$ be distinct. Write $t_{uv}$ for the number of $4$-cycles containing $c_0$ with at least one vertex in each of $J_u\setminus c_0$ and $J_v\setminus c_0$, and $d_u = \deg(u)$, $d_v=\deg(v)$.
        \begin{itemize}
        \item If $uv \notin E(G)$ then $t_{uv} \geq k^2 - k(d_u+d_v+2)$.
        \item If $uv\in E(G)$ then $t_{uv} \leq k^2 - k(d_u+d_v +2) -k + 2n^2 +3n$.
        \end{itemize}
    \end{claim*}

    \begin{poc}
    Let $c_1 \in J_u\setminus c_0$ and $c_3 \in J_v \setminus c_0$. There is a $4$-cycle $(c_0,c_1,c_2,c_3)$ in $\cC_k(G)$ precisely when there is some colouring $c_2$ that differs from $c_0$ at $u$ and $v$ only, and satisfies $c_2(v) = c_1(v)$ and $c_2(u) = c_3(u)$. An example of such a $4$-cycle is given in \cref{fig:c4p3}.
    
    Suppose that $uv \notin E(G)$.
    Since $c_3(u)$ must be distinct from both $c_0(u)$ and the colours of any neighbour of $u$ in $c_0$, there are $k - d_u - 1$ choices for $c_3$. Similarly, there are $k - d_v - 1$ choices for $c_1$, and hence the number of $4$-cycles containing $c_0$ and one vertex from each of $J_u$ and $J_v$ is
    \[
    k^2 - k(d_u + d_v + 2) + d_ud_v + d_u + d_v + 1 \geq k^2 - k(d_u + d_v + 2).
    \]
    Next suppose that $uv\in  E(G)$, and suppose first that we choose $c_3(u)$ to be a colour not used by $c_0$. Then there are $(k - n)$ choices for $c_3$, and since $c_2(v)$ must be distinct from each of $c_0(v)$, $c_3(u)$ and the colours of each neighbour of $v$ in $c_0$, this leaves $(k - d_v 
     - 2)$ choices for $c_1$ for a total of $(k - n)(k - d_v 
     - 2)$ pairs $(c_1, c_3)$. Similarly, if we choose $c_1$ first, we count $(k - n)(k - d_u 
     - 2)$ colour pairs. Any colour pair in which both $c_1(v)$ and $c_3(u)$ are selected from the set of $k-n$ colours not used by $c_0$ is counted twice above, so the total number of colour pairs in which at least one of $c_1(v)$ and $c_3(u)$ is a colour not used by $c_0$ is 
    \begin{align*}
        &\quad\ (k-n)[(k-d_u -2) + (k-d_v-2)- (k-n-1)] \\
        &= k^2 - k(d_u + d_v + 3) + n(d_u + d_v - n +3)\\
        &\leq k^2 -k(d_u + d_v + 3) +n^2 +3n.
    \end{align*}
    Since there are at most $n^2$ ways to choose $c_1(v)$ and $c_3(u)$ from colours used by $c_0$, we have the desired bound on $t_{uv}$.
    \end{poc}

     We now build a candidate graph $G_c$ from each vertex $c$ of $\cC_k(G)$ by considering pairs of cliques $J_u, J_v \in \cJ(c)$ and adding the edge $uv$ in $E(G_c)$ whenever $t_{uv} \geq k^2 - k(k-|J_u|+k-|J_v| + 2)$. 
     When $k> 5n^2$ we have $2n^2+3n < k$ for all positive $n$, so when $c$ is rainbow, $t_{uv} \geq k^2 - k(d_u+d_v+2)$ if and only if $uv \notin E(G)$. Furthermore, when $c$ is rainbow, we have $d_u = k - |J_u|$ for each $u \in V(G)$. Substituting this term into our formula for $t_{uv}$, we see that our candidate graph $G_c$ is isomorphic to $G$ whenever $c$ is rainbow. 
     Since \cref{lem:3n2_bound} guarantees that the majority of vertices in $\cC_k(G)$ are rainbow, more than half of these candidates are isomorphic to $G$, and so $G$ can be reconstructed by majority vote.
\end{proof}

\section{Finite families of polynomials}\label{sec:finite}

We now work towards proving \cref{thm:no_finite_family}. First, we show that for any finite collection $\cF$ of connected graphs, the polynomials $\{\pi_G^{(H)}(k)\colon H\in \cF\}$ cannot distinguish all graphs. This is implied by the following result.

\begin{lemma}\label{lem:connected_H_not_invariant}
    For each natural number $m$, there is a pair of non-isomorphic graphs $G$, $G'$ such that for every connected graph $H$ with at most $m$ edges, $\pi_G^{(H)} = \pi_{G'}^{(H)}$.
\end{lemma}

\begin{proof}
    Let $m$ be given, and choose any natural number $n > m+1$.
    Consider the graph $G_0$ obtained from the path $v_1,\dotsc, v_{3n}$ by adding a new vertex $v$ adjacent to $v_{n-1}$ and $v_n$, and a new vertex $v'$ adjacent to $v_n$ and $v_{n+1}$.
    Define $G$ and $G'$ to be the subgraphs of $G_0$ induced by $\{ v,v_1,\dotsc, v_{3n}\}$ and $\{v', v_1,\dotsc,v_{3n}\}$, respectively (see \cref{fig:G_and_Gprime}).

    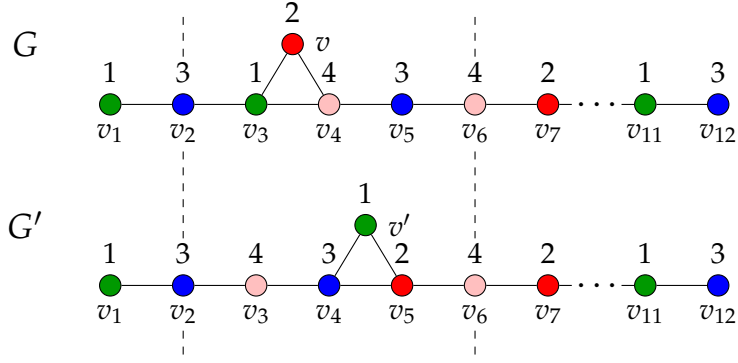
\begin{figure}[ht]
     \centering
     \begin{tikzpicture}[
    mycircle/.style={
        circle,
        draw=black,
        fill=black,
        fill opacity = 1,
        inner sep=0pt,
        minimum size=8pt,
        font=\small},
    nocircle/.style={
        circle,
        draw=black,
        fill=white,
        fill opacity = 0,
        inner sep=0pt,
        minimum size=0pt,
        font=\small},
    myarrow/.style={-},
    dottedarrow/.style={-,dashed},
    thiccarrow/.style={-,line width=0.9pt},
    node distance=1.2cm and 1.5cm
]

\tikzstyle{every node}=[font=\small]

\begin{scope}[scale=0.8]
\begin{scope}
    \node[mycircle, fill=green!60!black, label={$1$}, label={[yshift=-0.8cm]$v_1$}] (1) at (0,0){};
    \node[mycircle, fill=blue, label={$3$}, label={[yshift=-0.8cm]$v_2$}] (2) at (1.2,0){};
    \node[mycircle, fill=green!60!black, label={$1$}, label={[yshift=-0.8cm]$v_3$}] (3) at (2.4,0){};
    \node[mycircle, fill=pink, label={$4$}, label={[yshift=-0.8cm]$v_4$}] (4) at (3.6,0){};
    \node[mycircle, fill=red, label={$2$}, label={[]right:$v$}] (v) at (3,1){};
    \node[mycircle, fill=blue, label={$3$}, label={[yshift=-0.8cm]$v_5$}] (5) at (4.8,0){};
    \node[mycircle, fill=pink, label={$4$}, label={[yshift=-0.8cm]$v_6$}] (6) at (6,0){};
    \node[mycircle, fill=red, label={$2$}, label={[yshift=-0.8cm]$v_7$}] (7) at (7.2,0){};

    \node[nocircle] (ghost1) at (7.6,0){};
    \node[nocircle] (ghost2) at (8.4,0){};

    \node[mycircle, fill=green!60!black, label={$1$}, label={[yshift=-0.8cm]$v_{11}$}] (11) at (8.8,0){};
    \node[mycircle, fill=blue, label={$3$}, label={[yshift=-0.8cm]$v_{12}$}] (12) at (10,0){};

    \path[every node/.style={font=\sffamily\small}]
        (1) edge [color=black] (2)
        (2) edge [color=black] (3)
        (3) edge [color=black] (4)
        (3) edge [color=black] (v)
        (v) edge [color=black] (4)
        (4) edge [color=black] (5)
        (5) edge [color=black] (6)
        (6) edge [color=black] (7)
        (7) edge [color=black] (ghost1)
        (ghost2) edge [color=black] (11)
        (11) edge [color=black] (12);
    
    \node[scale=1.2] () at (8.05,0) {$\dotsc$};

    \node[scale=1.2] () at (-1.4,1) {$G$};
\end{scope}

\begin{scope}[yshift=-3cm]
    \node[mycircle, fill=green!60!black, label={$1$}, label={[yshift=-0.8cm]$v_1$}] (1) at (0,0){};
    \node[mycircle, fill=blue, label={$3$}, label={[yshift=-0.8cm]$v_2$}] (2) at (1.2,0){};
    \node[mycircle, fill=pink, label={$4$}, label={[yshift=-0.8cm]$v_3$}] (3) at (2.4,0){};
    \node[mycircle, fill=blue, label={$3$}, label={[yshift=-0.8cm]$v_4$}] (4) at (3.6,0){};
    \node[mycircle, fill=green!60!black, label={$1$}, label={[]right:$v'$}] (v') at (4.2,1){};
    \node[mycircle, fill=red, label={$2$}, label={[yshift=-0.8cm]$v_5$}] (5) at (4.8,0){};
    \node[mycircle, fill=pink, label={$4$}, label={[yshift=-0.8cm]$v_6$}] (6) at (6,0){};
    \node[mycircle, fill=red, label={$2$}, label={[yshift=-0.8cm]$v_7$}] (7) at (7.2,0){};

    \node[nocircle] (ghost1) at (7.6,0){};
    \node[nocircle] (ghost2) at (8.4,0){};

    \node[mycircle, fill=green!60!black, label={$1$}, label={[yshift=-0.8cm]$v_{11}$}] (11) at (8.8,0){};
    \node[mycircle, fill=blue, label={$3$}, label={[yshift=-0.8cm]$v_{12}$}] (12) at (10,0){};

    \path[every node/.style={font=\sffamily\small}]
        (1) edge [color=black] (2)
        (2) edge [color=black] (3)
        (3) edge [color=black] (4)
        (4) edge [color=black] (v')
        (v') edge [color=black] (5)
        (4) edge [color=black] (5)
        (5) edge [color=black] (6)
        (6) edge [color=black] (7)
        (7) edge [color=black] (ghost1)
        (ghost2) edge [color=black] (11)
        (11) edge [color=black] (12);
    
    \node[scale=1.2] () at (8.05,0) {$\dotsc$};

    \node[scale=1.2] () at (-1.4,1) {$G'$};
\end{scope}

\draw [dashed] (1.2,1.45) -- (1.2, 0.66);
\draw [dashed] (1.2, -0.78) -- (1.2, -2.22);
\draw [dashed] (1.2, -3.75) -- (1.2,-4.2);

\draw [dashed] (6,1.45) -- (6, 0.66);
\draw [dashed] (6, -0.78) -- (6, -2.22);
\draw [dashed] (6, -3.75) -- (6,-4.2);

\end{scope}
\end{tikzpicture}
     \vspace{-1cm}
     \caption{A $k$-colouring $c$ of $G$ for $n=4$, and the corresponding $k$-colouring $f_X(c)$ of $G'$. In this example $v_{n-t_X} = v_2$ and $v_{n+t_X} = v_6$, so $f$ recolours only vertices in the segment between $v_2$ and $v_6$.}
     \label{fig:G_and_Gprime}
 \end{figure}
    
    Fix any connected graph $H$ with at most $m$ edges, and some natural number $k$. We will show that $\pi_G^{(H)}(k) = \pi_{G'}^{(H)}(k)$ for each $k$ by finding a bijection $f$ between induced copies of $H$ in $\cC_k(G)$ and in $\cC_k(G')$. For each copy $X$ of $H$ in $\cC_k(G)$, we will define $f(X)$ in terms of an intermediate map $f_X$ that transforms colourings of $G$ into colourings of~$G'$.

    We first construct $f_X$. For $G^\ast\in \{G,G'\}$, let $X^\ast$ be an induced copy of $H$ in $\cC_k(G^\ast)$. Say that a vertex $v\in V(G^\ast)$ corresponds to an edge $c_1c_2$ of $X^\ast$ if it is the unique vertex at which the colourings $c_1$ and $c_2$ differ. Then, define $t_{X^\ast}$ to be the smallest positive integer such that neither $v_{n-t_{X^\ast}}$ nor $v_{n+t_{X^\ast}}$ correspond to any edge of $X^\ast$. Such a value $t_{X^\ast} \in [n]$ exists because $X^\ast$ only has $m < n-1$ edges.
    
    Working in $G$, let $X$ be an induced copy of $H$ in $\cC_k(G)$. Fix $c\in V(X)$ and consider colours modulo $k$, taking $[k]$ as the set of representatives. We define $f_X(c)$ to colour each vertex $u$ of $G'$ by
    \begin{equation*}
        (f_X(c))(u) = \begin{cases}
            c(v_{n-t_X}) + c(v_{n+t_X}) - c(v_{n-i}) & \text{if $u = v_{n+i}$ for $i \in (-t_X, t_X)$,} \\
            c(v_{n-t_X}) + c(v_{n+t_X}) - c(v) & \text{if $u = v'$,}\\
            c(u) & \text{otherwise.}\\
        \end{cases}
    \end{equation*}
    \cref{fig:G_and_Gprime} provides an example of a colouring $c$ of $G$ and a corresponding colouring $f_X(c)$ of $G'$. 
    
    Next, let $X'$ be an induced copy of $H$ in $\cC_k(G')$, and fix $c' \in V(X')$. We similarly define $g_{X'}(c')$ to map colourings of $G'$ to colourings of $G$ by  
    \begin{equation*}
        (g_{X'}(c'))(u) = \begin{cases}
            c'(v_{n-t_{X'}}) + c'(v_{n+t_{X'}}) - c'(v_{n-i}) & \text{if $u = v_{n+i}$ for $i \in (-t_{X'}, t_{X'})$,} \\
            c'(v_{n-t_{X'}}) + c'(v_{n+t_{X'}}) - c'(v') & \text{if $u = v$,}\\
            c'(u) & \text{otherwise.}\\
        \end{cases}
    \end{equation*}
    Let $f(X)$ be the subgraph of $\cC_k(G')$ induced by $f_X(V(X))$, and let $g(X')$ be the subgraph of $\cC_k(G)$ induced by $g_{X'}(V(X'))$.

    For each $c \in V(X)$, we observe that $f_X(c)$ is indeed a proper colouring of~$G'$.  We now show that $X$ and $f(X)$ are isomorphic. 
Since $H$ is connected and neither $v_{n-t_X}$ nor $v_{n+t_X}$ correspond to an edge of $X$, the colours of $v_{n-t_X}$ and $v_{n+t_X}$ are constant across all colourings in $V(X)$.
Hence, the value $c(v_{n-t_X}) + c(v_{n+t_X})$ is the same across each colouring $c$ in $X$. 
It is then straightforward to check that two colourings $c_1$ and $c_2$ are adjacent in $X$ if and only if $f_X(c_1)$ and $f_X(c_2)$ are adjacent in~$f(X)$. This implies that $t_X = t_{f(X)}$, and that $f(X)$ is an induced copy of $H$ in $\cC_k(G')$. Making symmetric observations about $g_{X'}$, we now observe that $g_{f(X)}$ is the inverse of $f_X$, and that therefore $f$ is also invertible with inverse $g$. That is, $f$ is a bijection between induced copies of $H$ in $\cC_k(G)$ and induced copies of $H$ in $\cC_k(G')$, and the result follows.
\end{proof}

We now extend this result to finite collections of disconnected graphs via a standard argument.

\begin{lemma}\label{lem:discon}
Let $H$ be a graph with connected components $R_1,\dotsc, R_t$.
Then there is a finite collection $\cF$ of connected graphs such that $\{ \pi_G^{(F)} \colon F\in \cF \}$ uniquely determines $\pi_G^{(H)}$ for every graph $G$.
\end{lemma}
\begin{proof}
We proceed by induction on the number of connected components $t$ of $H$, with the base case being when $H$ is any connected graph.
Suppose $H$ has at least $t>1$ components $R_1,\dotsc, R_t$.

The product $\pi_G^{(R_1)}(k)\dotsb \pi_G^{(R_t)}(k)$ counts the number of tuples $(\rho_1,\dotsc,\rho_t)$ of injective maps $\rho_i\colon V(R_i)\rightarrow V(\cC_k(G))$ such that for each $i \in [t]$, $\rho_i(V(R_i))$ induces a copy of $R_i$ in $\cC_k(G)$.
Fix such a tuple of injective maps and let $F$ be the subgraph of $\cC_k(G)$ induced by the images of the maps, i.e. by the vertices in $\bigcup_{i} \rho_i(V(R_i))$.
Notice that for a fixed graph $H$, there are finitely many possible isomorphism classes for the graph $F$ (all such graphs have at most $|H|$ vertices), and that $F$ is either isomorphic to $H$, or else has fewer connected components than $H$.
If we fix such an isomorphism class $F$, its contribution to the count $\pi_G^{(R_1)}(k)\dotsb \pi_G^{(R_t)}(k)$ is precisely~$\pi_G^{(F)}(k)$ times the number $N(F,H)$ of tuples $(\rho_1,\dotsc, \rho_t)$ which produce the vertices of the same copy of~$F$ in $G$.
Hence, letting $\cF^*$ be the family of non-isomorphic graphs other than $H$ that can be obtained in this way, the following equality holds:
\begin{equation}
\label{eq:induct}
\pi^{(H)}_G = \pi^{(R_1)}_G \dotsb \pi^{(R_t)}_G -\sum_{F \in \cF^*} N(F,H)\cdot \pi^{(F)}_G. 
\end{equation}
By the induction hypothesis, each graph $F\in \cF^*$ has a finite collection of connected graphs $\cF_F$ (possibly $\cF_F=\{F\}$) which determine $\pi_G^{(F)}$, and so the finite family $\cF = \bigcup_{F\in \cF^*} \cF_F$ of connected graphs determines $\pi_G^{(H)}$.
\end{proof}

\begin{remark}\label{rem:finitepoly}
For connected $H$, the preceding proof can still be used to find a finite family of graphs $\cF(H)$, not containing $H$ or depending on $G$, such that the collection $\{ \pi_G^{(F)} \colon F\in \cF(H)\}$ determines $\pi_G^{(H)}$ (when $H$ is disconnected this comes directly from the proof with formula given by \cref{eq:induct}). Namely, run the proof with a disconnected graph $H^+$ that contains $H$ as one connected component. Then, isolating the term $\pi_G^{(H)}$ (which now occurs among the components and possibly in $\cF^\ast$) in \cref{eq:induct} gives the relevant formula.
\end{remark} 

\finitefamily*

\begin{proof}
Let $H_1,\dotsc,H_t$ be a finite family of graphs.
By \cref{lem:discon}, there is a finite collection $\cF$ of connected graphs such that for every graph $G$ the generalised chromatic polynomials $\pi_G^{(H_1)}, \dotsc, \pi_G^{(H_t)}$ only depend on $\{ \pi_G^{(F)} \colon F\in \cF\}$.
The theorem now follows from choosing $m$ in \cref{lem:connected_H_not_invariant} to be larger than $\max \{ |E(F)| \colon F\in \cF\}$.
\end{proof}

\section{Open problems}
Asgarli et al. conjectured that $\pi_G^{(K_2)}$ determines the chromatic polynomial $\pi_G^{(K_1)}$, that is, if $\pi_{G_1}^{(K_2)} = \pi_{G_2}^{(K_2)}$ then $\pi_{G_1}^{(K_1)} = \pi_{G_2}^{(K_1)}$ \cite[Conjecture 5.2]{AKLR24}. This remains unverified, but in light of \cref{rem:finitepoly} we ask a broader question.

\begin{problem}\label{pb:pair}
    For which graphs $H$ does there exist an $H'$ such that $\pi_G^{(H')}$ determines~$\pi_G^{(H)}$ for every graph $G$?
\end{problem}
Define the \defn{graph product} of $G_1, G_2$ to be the graph $G_1 \sq G_2$ on vertex set $V(G_1) \times V(G_2)$ with adjacencies between vertices $(a_1,b_1)$, $(a_2,b_2)$ if and only if $a_1 =a_2$ and $b_1b_2\in E(G)$ or $b_1=b_2$ and $a_1a_2\in E(G)$.
Notice that for any graph $G$ on $n$ vertices and any integer $k$, the graph $\cC_k(G)$ is obtained from the product $Q(k,n)\coloneqq K_k \sq \dotsb \sq K_k$ of $n$ copies of $K_k$ by removing vertices corresponding to $k$-colourings of vertices of $G$ which are not proper in $G$.
As such, every graph $H$ whose polynomial $\pi_G^{(H)}$ is nonzero for some graph $G$ must be an induced subgraph of $Q(k,n)$ for some $k$ and $n$.
\cref{pb:pair} is therefore only interesting when $H$ is an induced subgraph of $Q(k,n)$.

It is well known that the chromatic polynomial does not distinguish all graphs, and 
we have shown in \cref{thm:no_finite_family} that no finite family of generalised chromatic polynomials suffices to distinguish all graphs.  But what about typical graphs?
Bollob\'as, Pebody and Riordan~\cite[Conjecture 2]{BOLLOBAS20} raised the intriging conjecture that almost every graph is determined by its chromatic polynomial; they asked the same question for the Tutte polynomial (which is a stronger invariant).
In a similar vein, we ask a weakening of their chromatic polynomial conjecture for finite families of generalised chromatic polynomials.
Let $\cG(n,p)$ be the random graph on $n$ vertices obtained by sampling each edge independently with probability $p$.
\begin{problem}
    Is there a finite set of graphs $H_1,\dotsc, H_t$ such that, for almost every \hbox{$G\in \cG(n,\frac12)$}, if  $\pi_{G'}^{(H_i)} = \pi_{G}^{(H_i)}$ for all $i$ then $G'\cong G$?
\end{problem}

\bibliography{chromaticpolys}
\bibliographystyle{style_edited}
\end{document}